\documentclass{amsart}
\usepackage{amssymb,amsmath,amsthm,graphics,amscd,amsfonts}

\theoremstyle{plain}
\newtheorem{theorem}{Theorem}[section]
\newtheorem{proposition}[theorem]{Proposition}

\newtheorem{lemma}[theorem]{Lemma}

\newtheorem{definition}[theorem]{Definition}

\newcommand{\bfC}{{\mathbb C}}

\newcommand{\bari}{{\overline i}}
\newcommand{\barj}{{\overline j}}
\newcommand{\bark}{{\overline k}}
\newcommand{\barl}{{\overline \ell}}

\newcommand{\baru}{{\overline u}}

\newcommand{\barq}{{\overline q}}
\newcommand{\barz}{{\overline z}}
\newcommand{\barv}{{\overline v}}

\newcommand{\mapright}[1]{\smash{\mathop{   \hbox to 0.7cm{\rightarrowfill}}
  \limits^{#1}}}

\newcommand{\Ker}{{\rm Ker}}

\def\a{\alpha}

\def\bp{\overline{\partial}}

\def\grad{\mathrm{grad}}

\def\J{\mathcal J}

\def\p{\partial}
\def\bp{\overline{\partial}}

\begin{document}
\title{Holomorphic vector fields and 
perturbed extremal K\"ahler metrics}
\author{Akito Futaki}
\address{Department of Mathematics, Tokyo Institute of Technology, 2-12-1,
O-okayama, Meguro, Tokyo 152-8551, Japan}
\email{futaki@math.titech.ac.jp}

\date{January 15, 2007 }

\begin{abstract} 
We prove a theorem which
asserts that the Lie algebra of all holomorphic vector fields on a compact 
K\"ahler manifold with a perturbed extremal metric has the structure 
similar to the case of an unperturbed extremal K\"ahler metric proved by Calabi.
 \end{abstract}



\maketitle

\section{Introduction}
Let $M$ be a compact symplectic manifold with symplectic form $\omega$.
On the space $\J$ of all $\omega$-compatible complex structures $J$ there is 
a natural symplectic form with respect to which the scalar curvature $S(J)$ of the
K\"ahler manifold $(M, \omega, J)$ becomes a moment map for the action of
the group of all Hamiltonian diffeomorphisms of $(M, \omega)$ acting on $\J$
(c.f. \cite{donaldson97}, \cite{fujiki92}). This means that the problem of finding
extremal K\"ahler metrics can be set in the framework of stability in the sense of
geometric invariant theory. 
It was shown in \cite{futaki06} that,  perturbing the symplectic form on $\J$ and 
the scalar curvature incorporating with the higher Chern classes and with a small real 
parameter $t$, the perturbed scalar curvature $S(J,t)$
becomes  a moment map with respect to the perturbed symplectic form on $\J$.
Note that the unperturbed scalar curvature is the trace
of the first Chern class. See section 2 for the precise definitions. 

Recall that a K\"ahler metric $g$ is called an extremal K\"ahler metric if 
the $(1,0)$-part of the gradient vector field of the scalar curvature $S$
$$ \mathrm{grad}' S = g^{i\barj}\frac{\p S}{\p \barz^j}\frac{\p}{\p z^i}$$
is a holomorphic vector field. 
Extremal K\"ahler metrics are critical points of two functionals. One is the so-called
the Calabi functional. This is a functional $\Psi$ on the space $\mathcal K_{\omega_0}$ 
of all K\"ahler forms in a
fixed de Rham class $\omega_0$ with fixed complex structure $J$. 
If $\omega \in \mathcal K_{\omega_0}$ and $S(\omega)$ denotes
the scalar curvature of $\omega$ then 
$$ \Psi(\omega) = \int_M S(\omega)^2 \omega^m$$
where $m = \dim_{\bfC}M$. 
Calabi originally defined extremal K\"ahler metrics to be the critical points of $\Psi$. 
The other functional $\Phi$ is defined on $\J$. If $S(J)$ denotes the scalar curvature
of the K\"ahler manifold $(M,\omega,J)$ for $J \in \J$ then
$$ \Phi(J) = \int_M S(J)^2\omega^m.$$
It is easy to see that the extremal K\"ahler metrics are exactly the critical points of $\Phi$ 
from the fact that the scalar curvature is the moment
map on $\J$ for the action of Hamiltonian diffeomorphisms as mentioned above.

Inspired by a work of Bando \cite{bando83} the author defined in \cite{futaki06} 
perturbed extremal
K\"ahler metrics as follows: the K\"ahler metric $g$ for $(M, \omega, J)$ is called a
perturbed extremal K\"ahler metric if 
the $(1,0)$-part of the gradient vector field
$$ \mathrm{grad}' S(J,t) = g^{i\barj}\frac{\p S(J,t)}{\p \barz^j}\frac{\p}{\p z^i}$$
is a holomorphic vector field. From the fact that $S(J,t)$ becomes a moment map
on $\J$ with respect to the perturbed symplectic structure, one can see that
the critical points of the functional
$$\Phi(J) = \int_M S(J,t)^2 \omega^m$$
are $J$'s for which the K\"ahler metric of $(M,\omega,J)$ is a perturbed extremal K\"ahler
metric. 
However it is not true for $t \ne 0$ that perturbed extremal K\"ahler metrics are the critical points
of the functional $\Psi$ on $\mathcal K_{\omega_0}$ defined by
$$ \Psi(\omega) = \int_M S(\omega,t)^2 \omega^m$$
where $S(\omega,t)$ is the perturbed scalar curvature of $(M, \omega, J)$, 
see Remark 3.3 in \cite{futaki06}. This is the significant difference between the perturbed case
and the unperturbed case.

In \cite{xwang04} Xiaowei Wang explains how one gets the decomposition theorem of Calabi 
\cite{calabi85} for
the structure of the Lie algebra of all holomorphic vector fields on compact K\"ahler
manifolds with extremal K\"ahler metrics in the finite dimensional setting of the 
framework of the moment maps, see also \cite{futaki05}. 
On the other hand Lijing Wang \cite{Lijing06} explains
how one gets the Hessian formulae for the Calabi functional and the functional $\Phi$ in the finite dimensional
setting of the framework of moment maps. Recall that the Hessian formula for the Calabi
functional plays the key role for the proof of Calabi's decomposition theorem of the Lie
algebra of all holomorphic vector fields on compact K\"ahler manifolds with extremal
K\"ahler metrics. Because of the above mentioned difference between the perturbed case
and the unperturbed case, one can not expect that the same proof as the unperturbed case
by Calabi can be applied to the perturbed case. 
The purpose of this paper is to see L.-J. Wang's finite dimensional arguments provide us a
rigorous proof of Calabi's decomposition theorem for compact K\"ahler manifolds
with perturbed extremal K\"ahler metrics. Thus we obtain a similar statement of the 
decomposition theorem:

\begin{theorem}\label{main1} Let $M$ be a compact K\"ahler manifold with a perturbed
extremal K\"ahler metric. Let $\mathfrak h(M)$ be the Lie algebra of all holomorphic vector
fields and $\mathfrak k$ be the real Lie algebra of all Killing vector fields of $M$.
Then
\begin{enumerate}
\item[(a)] \ \ $\mathfrak h_0(M) := \mathfrak k \otimes \bfC$ is the maximal reductive subalgebra 
of $\mathfrak h(M)$.
\item[(b)]\ \ The $(1,0)$-part of the gradient vector field
$$ \mathrm{grad}' S(J,t) = g^{i\barj}\frac{\p S(J,t)}{\p \barz^j}\frac{\p}{\p z^i}$$
of $S(J,t)$ belongs to the center of $\mathfrak h_0(M)$. 
\item[(c)]\ \ $\mathfrak h(M)$ has the structure of semi-direct decomposition
$$\mathfrak h(M) = \mathfrak h_0(M) + \sum_{\lambda\ne 0} \mathfrak h_{\lambda}(M)$$
where $\mathfrak h_{\lambda}(M)$ is the $\lambda$-eigenspace of the adjoint action
of $\mathrm{grad}'S(J,t) $.
\end{enumerate}
\end{theorem}

We will follow the arguments of L.-J. Wang almost word for word.

Throughout this paper Hermitian inner products are anti-linear in the first component and linear in
the second component. 

\section{Perturbed extremal K\"ahler metric}

Let $M$ be a compact symplectic manifold of dimension $2m$
with symplectic form $\omega$, $\J$ the space of all $\omega$-compatible
complex structures on $M$. Then for each $J \in \J$, $(M,\omega, J)$ becomes a
K\"ahler manifold. For a pair $(J, t)$, $t$ being a small real number, we define a
smooth function $S(J,t)$ on $M$ by
\begin{equation}\label{S(J,t)}
 S(J,t)\, \omega^m = c_1(J) \wedge \omega^{m-1} + 
t c_2(J) \wedge \omega^{m-2} + \cdots + t^{m-1} c_m(J) 
\end{equation}
where $c_i(J)$ is the $i$-the Chern form defined by
\begin{equation}\label{Chern}
\det(I + \frac{i}{2\pi}t\Theta) = 1 + tc_1(J) + \cdots + t^m c_m(J),
\end{equation}
$\Theta$ being the curvature form with respect to $\omega$.
Note that we use $S(J,t)$ in place of $S(J,T)/2m\pi$ in \cite{futaki06} to avoid clumsy constant
$1/2m\pi$. 
\begin{definition}\label{PExt} The K\"ahler metric $g$ of the K\"ahler manifold
$(M, J, \omega)$ is called a $t$-perturbed extremal K\"ahler metric or simply 
perturbed extremal metric if 
\begin{equation}\label{grad}
\grad'S(J,t) = \sum_{i,j = 1}^m g^{i\barj}\frac{\p S(J,t)}{\p \barz^j} \frac{\p}
{\p z^i}
\end{equation}
is a holomorphic vector field.
\end{definition}

The following was proved in \cite{futaki06}, Proposition 3.2.
\begin{proposition}\label{critical}
The critical points of the functional $\Phi$ on $\J$ defined by
\begin{equation}\label{Phi}
\Phi(J) = \int_M S(J,t)^2 \omega^m
\end{equation}
are the perturbed extremal K\"ahler metrics.
\end{proposition}
The proof of this proposition essentially follows from the fact that the
perturbed scalar curvature $S(J,t)$ gives the moment map for the action
of the group of Hamiltonian diffeomorphisms with respect to a perturbed
symplectic structure on $\J$. This perturbed symplectic structure is described as 
follows. 
The tangent space of $\J$ at $J$ is identified with a subspace of $\mathrm{Sym}(\otimes^2
T^{\prime\prime\ast}M)$. For a small real number $t$, we define an Hermitian structure on 
$\mathrm{Sym}(\otimes^2
T^{\prime\prime\ast}M)$ by
\begin{equation}\label{inner}
(\nu, \mu)_t = \int_M mc_m(\overline{\nu}_{jk}
\,\mu^i{}_\barl \frac {\sqrt{-1}}{2\pi}\,
dz^k \wedge d\overline{z^{\ell}}, \omega\otimes I + \frac {\sqrt{-1}}{2\pi}\,t\Theta,
\cdots, \ \omega\otimes I + \frac {\sqrt{-1}}{2\pi}\,t\Theta)
\end{equation}
for $\mu$ and $\nu$ in the tangent space $T_J\mathcal J$, 
where $c_m$ is the polarization of the determinant viewed as a $GL(m,\bfC)$-invariant 
polynomial, i.e. $c_m(A_1, \cdots, A_m)$ is the coefficient of $m!\, t_1 \cdots t_m$ in $\det(t_1A_1 + \cdots +
t_mA_m)$, where $I$ denotes the identity matrix and $\Theta = \bp(g^{-1}\partial g)$ 
is the curvature form of the
Levi-Civita connection, and where 
$u_{jk}\mu^i_{\bar l}$  should be understood as
the endomorphism of  $T'_JM$ which sends 
$\partial/\partial z^j$  to  $u_{jk}\mu^i_{\bar l}{\partial/\partial z^i}$.
When $t = 0$, (\ref{inner}) gives the usual $L^2$-inner product. 
The perturbed symplectic form $\Omega_{J,t}$
at $J \in \J$ is then given by
\begin{eqnarray}\label{symp}
&&\Omega_{J,t}(\nu,\mu) =
\Re (\nu,\sqrt{-1}\mu)_t \\
&&= \Re \int_M mc_m(\overline{\nu}_{jk}
\,\sqrt{-1}\mu^i{}_\barl \frac {\sqrt{-1}}{2\pi}\,
dz^k \wedge d\overline{z^{\ell}}, \omega\otimes I + \frac {\sqrt{-1}}{2\pi}\,t\Theta, \nonumber \\ 
&& \hspace{6.5cm}\cdots, \ \omega\otimes I + \frac {\sqrt{-1}}{2\pi}\,t\Theta) \nonumber
\end{eqnarray}
where $\Re$ means the real part. In \cite{futaki06} we proved the following:

\begin{theorem}[\cite{futaki06}] If $\delta J = \mu$ then
\begin{equation}\label{moment1}
\delta \int_M u\ S(J,t)\omega^m =  
\Omega_{J,t}(2\sqrt{-1}\nabla^{\prime\prime}\nabla^{\prime\prime}u,\mu).
\end{equation}
Namely the perturbed scalar curvature $S(J,t)$ gives a moment map with respect to the
perturbed symplectic form $\Omega_{J,t}$ for the action of the group of Hamiltonian
diffeomorphisms on $\J$.
\end{theorem}

Now we can prove Proposition \ref{critical}. From (\ref{moment1}) we have 
\begin{eqnarray}\label{moment1.5}
\delta \int_M S(J,t)^2\, \omega^m &=& 
2 \int_M S(J,t) \delta S(J,t)\, \omega^m \\
&=& 2\Omega_{J,t}(2\sqrt{-1}\nabla^{\prime\prime}\nabla^{\prime\prime}S(J,t), \mu).\nonumber
\end{eqnarray}
This shows that $J$ is a critical point if and only if 
\begin{equation}\label{moment2}
\nabla^{\prime\prime}\mathrm{grad}^{\prime}S(J,t) = 0,
\end{equation}
i.e. the K\"ahler metric of $(M, \omega, J)$ is a perturbed extremal K\"ahler
metric.

Let $\frak g$ be the complexification of the Lie algebra of the
group of Hamiltonian diffeomorphisms. Then $\frak g$ is simply the set of all
complex valued smooth functions $u$ with the normalization
$$ \int_M u\, \omega^m = 0$$
with the Lie algebra structure given by the Poisson bracket. 
The infinitesimal action of $u$ on $\J$ is given by $2i\nabla^{\prime\prime}\nabla^{\prime\prime}u$,
see Lemma 10 in \cite{donaldson97} or Lemma 2.3 in \cite{futaki06}.
Define $L : C^{\infty}(M)\otimes \bfC\ (\cong \frak g)\ \to C^{\infty}(M)\otimes \bfC $ by
\begin{eqnarray}\label{L}
&&(v,Lu)_{L^2} = (\nabla^{\prime\prime}\nabla^{\prime\prime}v,\nabla^{\prime\prime}
\nabla^{\prime\prime}u)_t \\
&& = \int_M mc_m(\barv_{jk}u^i{}_\barl \frac {\sqrt{-1}}{2\pi}\,
dz^k \wedge d\overline{z^{\ell}}, \omega\otimes I + \frac {\sqrt{-1}}{2\pi}\,t\Theta,
\cdots, \ \omega\otimes I + \frac {\sqrt{-1}}{2\pi}\,t\Theta). \nonumber
\end{eqnarray}
More explicitly $L$ is expressed as
\begin{equation}\label{L'}
Lu = mc_m(u^i{}_{\barl jk}\frac {\sqrt{-1}}{2\pi}\,
dz^k \wedge d\overline{z^{\ell}}, \omega\otimes I + \frac {\sqrt{-1}}{2\pi}\,t\Theta,
\cdots, \ \omega\otimes I + \frac {\sqrt{-1}}{2\pi}\,t\Theta)/\omega^m.
\end{equation}
We define ${\overline L} : C^{\infty}(M)\otimes \bfC \to C^{\infty}(M)\otimes \bfC $ by
$\overline{L}u := \overline{L \baru}$. Then $\overline{L}$ satisfies
\begin{eqnarray}\label{Lbar}
&&(v,\overline{L} u)_{L^2} = (\nabla^{\prime\prime}\nabla^{\prime\prime}\baru,
\nabla^{\prime\prime}\nabla^{\prime\prime}\barv)_t \\
&& = \int_M mc_m(u_{jk}\barv^i{}_\barl \frac {\sqrt{-1}}{2\pi}\,
dz^k \wedge d\overline{z^{\ell}}, \omega\otimes I + \frac {\sqrt{-1}}{2\pi}\,t\Theta,
\cdots, \ \omega\otimes I + \frac {\sqrt{-1}}{2\pi}\,t\Theta). \nonumber
\end{eqnarray}
and
\begin{equation}\label{Lbar'}
\overline{L} u =  mc_m(u_{jk}{}^i{}_{\barl}\frac {\sqrt{-1}}{2\pi}\,
dz^k \wedge d\overline{z^{\ell}}, \omega\otimes I + \frac {\sqrt{-1}}{2\pi}\,t\Theta,
\cdots, \ \omega\otimes I + \frac {\sqrt{-1}}{2\pi}\,t\Theta)/\omega^m.
\end{equation}

\begin{lemma}\label{deltaS}
If $v$ is a real smooth function and $\delta J = \nabla^{\prime\prime}\nabla^{\prime\prime}v$
then
$$ \delta S(J,t) = Lv + \overline{L}v.$$
\end{lemma}
\begin{proof} Let $u$ be also a real smooth function. Then by (\ref{moment1})
\begin{eqnarray*}
\int_M u\,\delta S(J,t) \omega^m &=& 
\Re (2\sqrt{-1}\nabla^{\prime\prime}\nabla^{\prime\prime}u,\sqrt{-1} \mu)_t \\
&=& (\nabla^{\prime\prime}\nabla^{\prime\prime}u, 
\nabla^{\prime\prime}\nabla^{\prime\prime}v)_t + 
(\nabla^{\prime\prime}\nabla^{\prime\prime}v,\nabla^{\prime\prime}\nabla^{\prime\prime}u)_t\\
&=& (u,Lv)_{L^2} + (u,\overline{L}v)_{L^2}.
\end{eqnarray*}
\end{proof}

\begin{lemma}\label{bracket}
Let $u$ and $v$ be real smooth functions and put 
$\mathcal X_u = 2\sqrt{-1}\nabla^{\prime\prime}\nabla^{\prime\prime}u$ and 
$\mathcal X_v = 2\sqrt{-1}\nabla^{\prime\prime}\nabla^{\prime\prime}v$. Then
we have
$$ \Omega_{J,t}(\mathcal X_u, \mathcal X_v) = (\{u,v\},S(J,t))_{L^2}.$$
\end{lemma}

\begin{proof}
Consider $\mathcal X_u$ and 
$\mathcal X_v$ 
as the infinitesimal action of real
Hamiltonian functions $u$ and $v$ on $\J$.
Since $S(J,t)$ gives an equivariant moment map
\begin{equation}\label{moment3}
 \int_M u S(\sigma J,t)\, \omega^m =  \int_M (\sigma^{-1\ast} u) S(J,t)\, \omega^m
 \end{equation}
for a Hamiltonian diffeomorphism  $\sigma$. If $\sigma$ is generated by the Hamiltonian
vector field of
a Hamiltonian function $v$ then (\ref{moment1}) and (\ref{moment3})
show
\begin{equation}\label{moment4}
\Omega_{J,t}(2\sqrt{-1}\nabla^{\prime\prime}\nabla^{\prime\prime}u, 2\sqrt{-1}\nabla^{\prime\prime}\nabla^{\prime\prime}v) = - \int_M S(J,t) \{v,u\}\omega^m.
\end{equation}
\end{proof}

\begin{lemma}\label{Lbar-L} For any smooth complex valued function $u$ we have
$$(\overline{L} - L) u = - \frac 12 (S(J,t)^{\a}u_{\a} - u^{\a}S(J,t)_{\a})$$
where $z^{\alpha}$'s are local holomorphic coordinates.
\end{lemma}
\begin{proof}It is sufficient prove when $u$ is a real valued function. Let $v$ be also a 
real valued smooth function. From (\ref{L}) and (\ref{Lbar}) we have
\begin{eqnarray*}
(v, \overline{L}u - Lu)_{L^2} &=& (\nabla^{\prime\prime}\nabla^{\prime\prime}u,
\nabla^{\prime\prime}\nabla^{\prime\prime}v)_t - 
(\nabla^{\prime\prime}\nabla^{\prime\prime}v,\nabla^{\prime\prime}\nabla^{\prime\prime}u)_t\\
&=& \overline{(\nabla^{\prime\prime}\nabla^{\prime\prime}v,
\nabla^{\prime\prime}\nabla^{\prime\prime}u)_t } - 
(\nabla^{\prime\prime}\nabla^{\prime\prime}v,
\nabla^{\prime\prime}\nabla^{\prime\prime}u)_t .
\end{eqnarray*}
It follows from this that
\begin{eqnarray*}
2\Re (\nabla^{\prime\prime}\nabla^{\prime\prime}v,
i\nabla^{\prime\prime}\nabla^{\prime\prime}u)_t
&=& i(\nabla^{\prime\prime}\nabla^{\prime\prime}v,
\nabla^{\prime\prime}\nabla^{\prime\prime}u)_t
+ \overline{i(\nabla^{\prime\prime}\nabla^{\prime\prime}v,
\nabla^{\prime\prime}\nabla^{\prime\prime}u)_t}\\
&=&-i(v,(\overline{L} - L)u)_{L^2}.
\end{eqnarray*}
Let $X_u$ denote the Hamiltonian vector field of $u$:
$ i(X_u)\omega = du$. Then $X_u = J\mathrm{grad}\, u$ and $\{u,S\} = X_uS$.
It then follows that
\begin{eqnarray*}
 (v,(\overline{L} - L)u)_{L^2} &=&
2i\Re(\nabla^{\prime\prime} \nabla^{\prime\prime} v, i\nabla^{\prime\prime} \nabla^{\prime\prime} u)_t \\
&=& \frac i2 \Re(\mathcal X_v, i\mathcal X_u)
= \frac i2\Omega_{J,t}(\mathcal X_v, \mathcal X_u)\\
&=& -\frac i2 (\{u,v\},S(J,t))_{L^2}
= \frac i2 (v,\{u, S(J,t)\})_{L^2}\\
&=& \frac i2 (v,X_uS(J,t))_{L^2} = \frac i2 \omega(v,g(X_u, J\mathrm{grad}\,S(J,t)))_{L^2}\\
&=& \frac i2 (v,du(J\mathrm{grad}\,S(J,t)))_{L^2}
= -\frac 12 (v,S(J,t)^{\alpha}u_{\alpha} - u^{\alpha} S(J,t)_{\alpha})_{L^2}.
  \end{eqnarray*}
\end{proof}
  
\begin{lemma}\label{d1}
Let $u$ be a real smooth function and suppose $\delta J = \nabla^{\prime\prime} \nabla^{\prime\prime}u$. Then
$$ \delta \int_M S(J,t)^2 \omega^m = 4(u,LS(J,t))_{L^2} = 4(u,\overline{L} S(J,t))_{L^2}.$$
\end{lemma}
\begin{proof} By (\ref{moment1.5})
\begin{eqnarray*}
\delta \int_M S(J,t)^2\omega^m &=& 2\Omega_{J,t}(2i\nabla^{\prime\prime} \nabla^{\prime\prime}
S(J,t), \nabla^{\prime\prime} \nabla^{\prime\prime}u) \\
&=& 4\Re (\nabla^{\prime\prime} \nabla^{\prime\prime} S(J,t), 
\nabla^{\prime\prime} \nabla^{\prime\prime}u)_t \\
&=& 2(\nabla^{\prime\prime} \nabla^{\prime\prime} S(J,t),
\nabla^{\prime\prime} \nabla^{\prime\prime}u)_t +
2(\nabla^{\prime\prime} \nabla^{\prime\prime}u, 
\nabla^{\prime\prime} \nabla^{\prime\prime} S(J,t))_t\\
&=& 2(u,LS(J,t))_{L^2} + 2(u,\overline{L}S(J,t))_{L^2}.
\end{eqnarray*}
But from Lemma \ref{Lbar-L} we have 
$$ \overline{L}S(J,t) = LS(J,t),$$
from which the lemma follows.
\end{proof}

\begin{lemma}\label{d2}
Suppose that $(\omega, J)$ is a perturbed extremal K\"ahler metric and thus that the gradient vector field of $S(J,t)$ is a holomorphic vector field. If $\delta J = \nabla^{\prime\prime} \nabla^{\prime\prime} u$ for a real smooth function $u$ then
$$ (\delta L)S(J,t) = -\frac12 L(S(J,t)^\alpha u_{\alpha} - u^{\alpha}S(J,t)_{\alpha}) = L(\overline{L} - L)u.$$
\end{lemma}

\begin{proof}\ \ Recall that by Lemma 2.3 in \cite{futaki06} 
$$ \mathcal L_XJ = 2i\nabla^{\prime\prime}_JX^{\prime} - 2i\nabla^{\prime}_JX^{\prime\prime} .$$
Therefore
\begin{eqnarray*}
\mathcal L_{JX}J &=& 2i\nabla^{\prime\prime}_JiX^{\prime} - 2i\nabla^{\prime}_J(-i)X^{\prime\prime}\\
&=& -2(\nabla^{\prime\prime}_JX^{\prime} - \nabla^{\prime}_JX^{\prime\prime}).\\
\end{eqnarray*}
This shows that $\mathcal L_{JX}J \in T_J \mathcal J$ corresponds to
$-2\nabla^{\prime\prime} \nabla^{\prime\prime} u \in 
\mathrm{Sym} \otimes^2T^{\prime\prime\ast}M$
via the identification $T_J \mathcal J \cong \mathrm{Sym} \otimes^2T^{\prime\prime\ast}M$.
Thus $\mathcal L_{-\frac12 JX_u}J$ corresponds to $\nabla^{\prime\prime}\nabla^{\prime\prime}u$.
On the other hand
\begin{equation}\label{d3}
\mathcal L_{\frac12 JX_u}\omega = d(i(\frac12 JX_u)\omega)
\end{equation}
and
\begin{eqnarray}\label{d4}
(i(\frac12 JX_u)\omega)(Y) &=& \omega(\frac12 JX_u,Y) = \omega(-\frac12 \mathrm{grad}u,Y)\\
&=& \omega(-\frac12 X_u,JY) = -\frac12 du\circ J = (d^cu)(Y) \nonumber
\end{eqnarray}
where $d^c = \frac i2(\bp - \p)$.
From (\ref{d3}) and (\ref{d4}) it follows that
\begin{equation}\label{d5}
\mathcal L_{\frac12 JX_u}\omega = dd^c u = i\p\bp u.
\end{equation}
Let $f_s$ is a flow generated by $-\frac 12 JX_u$.
Suppose that $S$ is a smooth function such that $\grad^{\prime}S$ is a holomorphic
vector field and that $S_s$ is a function such that
$$ \grad_s^{\prime} S_s = \grad^{\prime} S,\qquad \int_M S_s\, (f_{-s}^{\ast}\omega)^m = \int_M S\, \omega^m$$
where $\grad_s^{\prime} S_s$ is the $(1,0)$-part of the gradient vector field of $S_s$ with respect to 
$f_{-s}^{\ast}\omega$. It is easy to see that if $f_{-s}^{\ast} \omega = \omega + i\p\bp \varphi$
then $S_s = S + S^{\alpha} \varphi_{\alpha}$. Then (\ref{d5}) shows
\begin{equation}\label{d6}
S_s = S + sS^{\alpha}u_{\alpha} + O(s^2).
\end{equation} 
We have
\begin{equation}\label{d7}
L(f_sJ,\omega)f_s^{\ast}S_s = f_s^{\ast}(L(J,f_{-s}^{\ast}\omega)S_s) = 0.
\end{equation}
Takin the derivative of (\ref{d7}) with respect to t at $t=0$ we obtain
\begin{equation}\label{d8}
\delta L\cdot S + L(-\frac 12 (JX_u)S + S^{\alpha}u_{\alpha}) = 0.
\end{equation}
On the other hand
\begin{eqnarray}\label{d9}
JX_u\cdot S &=& g(JX_u,\grad\, S) = \omega(X_u,\grad\,S) = du(\grad\,S)\\
&=& (\p u + \bp u)(\nabla^{\prime} S + \nabla^{\prime\prime} S) = u^{\alpha}S_{\alpha}
+ S^{\alpha}u_{\alpha}.\nonumber
\end{eqnarray}
It follows from (\ref{d8}) and (\ref{d9}) that
\begin{eqnarray*}
\delta L\cdot S &=& - L(-\frac 12 (u^{\alpha}S_{\alpha} + S^{\alpha}u_{\alpha}) 
+ S^{\alpha}u_{\alpha})\\
&=& - L(\frac 12 (S^{\alpha}u_{\alpha} - u^{\alpha} S_{\alpha})).
\end{eqnarray*}
Applying this with $S = S(J,\omega)$ and using Lemma \ref{Lbar-L} complete the proof of Lemma \ref{d2}.
\end{proof}

\begin{theorem}\label{Hessian} Let $J$ be a critical point of $\Phi$, i.e. $(\omega, J)$
gives a perturbed extremal K\"ahler metric and $u$ be a real smooth function on $M$. 
Then the Hessian of $\Phi$ at $J$ in the
direction of $\nabla^{\prime\prime}\nabla^{\prime\prime} u$ and 
$\nabla^{\prime\prime}\nabla^{\prime\prime} v$is given by
$$ Hess(\Phi)_J(\nabla^{\prime\prime}\nabla^{\prime\prime} u, 
\nabla^{\prime\prime}\nabla^{\prime\prime} v) = 8(u,L\overline{L}v) =
8(u,\overline{L}Lv).$$
\end{theorem}

\begin{proof}\ \ Let $\delta J = \nabla^{\prime\prime}\nabla^{\prime\prime} v$.
By using Lemma \ref{d1}, Lemma \ref{d2} and Lemma \ref{deltaS}
successively one obtains
\begin{eqnarray*}
\mathrm{Hess}(\Phi)_J(\nabla^{\prime\prime}\nabla^{\prime\prime} u,
\nabla^{\prime\prime}\nabla^{\prime\prime} v)
&=& 4\delta (u,LS(J,t))\\
&=& 4(u,\delta L\cdot S(J,t) + L\delta S(J,t))\\
&=& 4(u,L(\overline{L} - L)v + L(L+ \overline{L})v)\\
&=& 8(u,L\overline{L}v).
\end{eqnarray*}
If one uses the third term in Lemma \ref{d1} and $\delta \overline{L} = L - \overline{L}$
then one gets the third term of Theorem \ref{Hessian}. This completes the proof.
\end{proof}

\section{Proof of Theorem \ref{main1}}

In this section we give a proof of Theorem \ref{main1}.
Suppose  that $g$ is a perturbed extremal K\"ahler metric on $(M, \omega, J)$.
Let $X$ be a holomorphic vector field and $\alpha$ be the dual $1$-form to $X$,
that is
$$ \alpha(Y) = g(X,Y), \qquad \alpha = \alpha_{\bari} d\overline{z}^i
= g_{j\bari}X^jd\barz^i. $$
Since $X$ is a holomorphic vector field
$$ \bp\alpha = (\nabla_\bari \alpha_{\barj} - \nabla_\barj \alpha_{\bari} )d\barz^i \wedge
d\barz^j = 0.$$
Let $\alpha = H\alpha + \bp\psi$ be the harmonic decomposition where $H\alpha$ denotes
the harmonic part. Then
\begin{eqnarray*}
L\psi &=& mc_m (\psi^i{}_{{\barl}jk}\frac {\sqrt{-1}}{2\pi}\,
dz^k \wedge d\overline{z^{\ell}}, \omega\otimes I + \frac {\sqrt{-1}}{2\pi}\,t\Theta,
\cdots, \ \omega\otimes I + \frac {\sqrt{-1}}{2\pi}\,t\Theta)\\
&=& mc_m ((X^i - (H\alpha)^i)_{\barl jk}\frac {\sqrt{-1}}{2\pi}\,
dz^k \wedge d\overline{z^{\ell}}, \omega\otimes I + \frac {\sqrt{-1}}{2\pi}\,t\Theta,\\
&&\hspace{7cm}\cdots, \ \omega\otimes I + \frac {\sqrt{-1}}{2\pi}\,t\Theta)\\
&=& - mc_m((H\alpha)^i{}_{\barl jk}\frac {\sqrt{-1}}{2\pi}\,
dz^k \wedge d\overline{z^{\ell}}, \omega\otimes I + \frac {\sqrt{-1}}{2\pi}\,t\Theta,
\cdots, \ \omega\otimes I + \frac {\sqrt{-1}}{2\pi}\,t\Theta)\\
&=& -mc_m((H\alpha)^i{}_{j\barl k} + (R_{j\barl}{}^i{}_p(H\alpha)^p)_k)\frac {\sqrt{-1}}{2\pi}\,
dz^k \wedge d\overline{z^{\ell}}, \omega\otimes I + \frac {\sqrt{-1}}{2\pi}\,t\Theta,\\
&&\hspace{7cm}
\cdots, \ \omega\otimes I + \frac {\sqrt{-1}}{2\pi}\,t\Theta).
\end{eqnarray*}
Note that being $\bp$-harmonic and being $\p$-harmonic are equivalent
on compact K\"ahler manifolds, and thus 
$$(H\alpha)_{\barq j} = \nabla_j(H\alpha)_{\barq} = 0.$$
This implies $(H\alpha)^i{}_j = 0$. It follows that
\begin{eqnarray}\label{pr1}
L\psi &=& -mc_m( R_{j\barl}{}^i{}_{p,k}(H\alpha)^p\frac {\sqrt{-1}}{2\pi}\,
dz^k \wedge d\overline{z^{\ell}}, \omega\otimes I + \frac {\sqrt{-1}}{2\pi}\,t\Theta,\\
&&\hspace{6cm}
\cdots, \ \omega\otimes I + \frac {\sqrt{-1}}{2\pi}\,t\Theta)\nonumber\\
&=& -mc_m( R_{j\barl}{}^i{}_{k,p}(H\alpha)^p\frac {\sqrt{-1}}{2\pi}\,
dz^k \wedge d\overline{z^{\ell}}, \omega\otimes I + \frac {\sqrt{-1}}{2\pi}\,t\Theta,\nonumber\\
&&\hspace{6cm}
\cdots, \ \omega\otimes I + \frac {\sqrt{-1}}{2\pi}\,t\Theta)\nonumber\\
&=& - (H\alpha)^p\nabla_pS(J,t) = - (H\alpha)_{\barq}\nabla^{\barq}S(J,t)\nonumber
\end{eqnarray}
where we have used the second Bianchi identity $R_{j\barl}{}^i{}_{p,k} = R_{j\barl}{}^i{}_{k,p}$
and 
\begin{eqnarray*}
\nabla_pS(J,t) &=& \nabla_p \frac 1t(c_m(\omega\otimes I + \frac i{2\pi}t\Theta) - \omega^m)\\
&=& \frac 1t \nabla_pc_m(\omega \otimes I + \frac i{2\pi}t\Theta)\\
&=& mc_m(R_{j\barl}{}^i{}_{k,p}\frac {\sqrt{-1}}{2\pi}\,
dz^k \wedge d\overline{z^{\ell}}, \omega\otimes I + \frac {\sqrt{-1}}{2\pi}\,t\Theta,\\
&&\hspace{4cm}
\cdots, \ \omega\otimes I + \frac {\sqrt{-1}}{2\pi}\,t\Theta).
\end{eqnarray*}

Note that $\nabla^{\barq}S(J,t)\frac{\p}{\p \barz^q}$ is a conjugate holomorphic
vector field and that $(H\alpha)_{\barq}d\barz^q$ is a conjugate holomorphic $1$-form
because $H\alpha$ is a $\p$-harmonic $(0,1)$-form. It follows from (\ref{pr1}) that
$L\psi = \mathrm{constant}$. But since $\int_M L\psi \omega^m = 0$ by (\ref{L})
we obtain $L\psi = 0$. This implies that $\grad^{\prime}\psi$ is a holomorphic
vector field. Then $(H\alpha)^i\frac{\p}{\p z^i} = X - \grad^{\prime}\psi$ is
also holomorphic. It then follows that
$$ \nabla_{\bark}(H\alpha)_{\barj} = 0.$$
But since $(H\alpha)$ is $\p$-harmonic we also have $\nabla_k(H\alpha)_{\barj} = 0$.
Thus $H\alpha$ is parallel.

This proves the direct sum decomposition as a vector space
$$ \frak h(M) = \frak a(M) +  \frak h^{\prime}(M)$$
where $$\frak h^{\prime}(M) = \{X \in \frak h(M)\ |\ X = \grad^{\prime}u\ \mathrm{for\ some}\
u \in C^{\infty}_{\bfC}(M)\}.$$
It is easy to see
$$ [\frak a(M), \frak a(M)] = 0;$$
$$ [\frak a(M),\frak h^{\prime}(M)] \subset \frak h^{\prime}(M);$$
$$ [\frak h^{\prime}(M),\frak h^{\prime}(M)] \subset \frak h^{\prime}(M).$$

Now by Theorem \ref{Hessian} we have $L\overline{L} = \overline{L}L$.
Thus $\overline{L}$ preserves $\Ker\, L$. Let $E_{\lambda}$ denote the
$\lambda$-eigenspace of $2\overline{L}|_{\Ker\,L}$. If $u \in E_{\lambda}$ then
$\grad^{\prime}u \in 
\frak h^{\prime}(M)$ and
\begin{eqnarray*}
 \lambda u &=& 2\overline{L} u \\
 &=& 2(\overline{L} - L)u \\
 &=& S(J,t)^{\alpha}u_{\alpha} - u^{\alpha}S(J,t)_{\alpha}.
 \end{eqnarray*}
 This implies $[\grad^{\prime}S(J,t),\grad^{\prime}u] = \lambda\, \grad^{\prime}u$.
 We put 
 $$\grad^{\prime}(E_{\lambda}) := \frak h_{\lambda}(M)\ \mathrm{for}\ \lambda \ne 0,$$
 $$\grad^{\prime}(E_0) := \frak h_0^{\prime}(M),$$
 $$\frak h_0 = \frak a(M) + \frak h_0^{\prime}(M).$$
 Then we obtain the decomposition 
 $$ \frak h(M) = \sum_{\lambda} \frak h_{\lambda}(M)$$
 where $\frak h_{\lambda}(M)$ is the $\lambda$-eigenspace of $\mathrm{ad} ( \grad^{\prime}
 S(J,t))$. Note that the real and imaginary parts of an element of $\frak a(M)$ are parallel and
 Killing and hence $[\grad^{\prime}S(J,t), \frak a(M)] = 0$.
 
 Finally since $E_0 = \Ker\,L \cap \Ker \overline{L}$, the real and imaginary parts
 are respectively in $E_0$, that is $E_0$ is the complexification of the purely imaginary
 functions $u$ such that $\grad^{\prime}u$ is holomorphic. The real parts of
 such $\grad^{\prime}u$'s are Killing vector fields, see Lemma 2.3.8 in \cite{futaki88}.
 The real parts of the elements of $\mathfrak a(M)$ are also Killing vector fields.
 Thus $\mathfrak h_0(M)$ is reductive.
 This completes the proof of Theorem \ref{main1}.

\end{document}